\overfullrule=0pt
\centerline {\bf A constrained approximation theorem for integral functionals on $L^p$}\par
\bigskip
\bigskip
\centerline {BIAGIO RICCERI}\par
\bigskip
\bigskip
{\bf Abstract.} Let $(T,{\cal F},\mu)$ be a $\sigma$-finite measure space, $E$ a separable real Banach space and $p\geq 1$. Given
a sequence of functions $f, f_1, f_2,...$ from $T\times E$ to ${\bf R}$, under general assumptions, we prove that, for each closed hyperplane
$V$ of $L^p(T,E)$, for each $u\in V$, and for each sequence $\{\lambda_n\}$ converging to $\int_Tf(t,u(t))d\mu$, there exists a sequence
$\{u_n\}$ in $V$ converging to $u$ and such that $\int_Tf_n(t,u_n(t))d\mu=\lambda_n$ for all $n$ large enough.\par
\bigskip
\bigskip
\bigskip
\bigskip
In the sequel, $(T,{\cal F},\mu)$ is a $\sigma$-finite non-atomic, complete measure space, $E$ is a separable real Banach space and $p\geq 1$.\par
\smallskip
As usual $L^p(T,E)$ denotes the space of all (equivalence classes of) measurable functions $u:T\to E$ such that
$\int_T\|u(t)\|_E^pd\mu<+\infty$, equipped with the norm
$$\|u\|_{L^p(T,E)}=\left(\int_T\|u(t)\|_E^pd\mu\right)^{1\over p}.$$
\smallskip
We also denote by ${\cal A}(T\times E)$ the class of all Carath\'eodory functions $g:T\times E\to {\bf R}$ such that, for each $u\in L^p(T,E)$,
the function $f(\cdot,u(\cdot))$ lies in $L^1(T)$. For each $g\in {\cal A}(T\times E)$, we set
$$J_g(u)=\int_Tg(t,u(t))d\mu$$
for all $u\in L^p(T,E)$.\par
\smallskip
Moreover, we denote by 
${\cal V}(L^p(T,E))$ the family of all sets $V \subseteq L^p(T,E)$ such that
$$V=\left\{u \in L^p(T,E): \Psi(u)=J_g(u)\right\}$$
where $\Psi \in (L^p(T,E))^*$, $g\in {\cal A}(T\times E)$ and $J_g$ is Lipschitzian on $L^p(T,E)$, with Lipschitz constant 
strictly smaller than $\|\Psi\|_{(L^p(T,E)*}$. Notice that, by a well-known representation theorem ([4], pp. 94-99), there exists a map
$P:T\to E^*$ such that the function $(t,x)\to P(t)(x)$ belongs to ${\cal A}(T\times E)$ and
$$\Psi(u)=\int_TP(t)(u(t))d\mu$$
for all $u\in L^p(T,E)$.
\smallskip
The aim of this paper is provide an answer to the following demanding approximation question:\par
\smallskip
Let $\{f_n\}$ be a sequence in ${\cal A}(T\times E)$ converging pointwise in $T\times E$ to a function $f\in {\cal A}(T\times E)$. When, for
each $V\in {\cal V}(L^p(T,E))$, for each $u\in V$ and for each sequence $\{\lambda_n\}$ in ${\bf R}$ converging to $J_{f}(u)$, there exists
a sequence $\{u_n\}$ in $V$ converging to $u$ and such that $J_{f_n}(u_n)=\lambda_n$ for all $n\in {\bf N}$ large enough ?
\smallskip
Our answer to this question is provided by the following\par
\medskip
THEOREM 1 . - {\it Let
let $\{f_n\}$ be a sequence of real-valued Carath\'eodory functions on $T\times E$ such that, for each $n\in {\bf N}$, one has
$$|f_n(t,x)|\leq M(t)+c\|x\|^p \eqno{(1)}$$
for a.e. $t\in T$ and for all $x\in E$, where $M\in L^1(T)$ and $c\geq 0$. Assume that $\{f_n\}$ converges pointwise in $T\times E$ to a 
Carath\'odory function $f$ such that
$$|f(t,x)|\leq M(t)\eqno{(2)}$$
for a.e. $t\in T$ and for all $x\in E$.
 Moreover, assume that the set
$$D=\{t\in T: f(t,\cdot)\hskip 5pt has\hskip 5pt  no\hskip 5pt global\hskip 5pt extrema\}$$
has a positive measure.\par
Then, for each  $V\in {\cal V}(L^p(T,E))$, for each $u\in V$ and for each sequence $\{\lambda_n\}$ converging
to $\int_Tf(t,u(t))d\mu$,
 there exist
a sequence $\{u_n\}$ in $V$ converging to $u$ and $\nu\in {\bf N}$ such that
$$\int_Tf_n(t,u_n(t))d\mu=\lambda_n$$
for all $n\geq \nu$.}
\smallskip
PROOF. Since $f_n$ is a Carath\'eodory function, in view of $(1)$,  we have $f_n\in {\cal A}(T\times E)$ and the
functional $J_{f_n}$ is continuous in $L^p(T,E)$, in view of the dominated convergence theorem. For each $u\in L^p(T,E)$, 
 since $$\lim_{n\to \infty}f_n(t,u(t))=f(t,u(t)),$$ the function
$f(\cdot,u(\cdot))$ is measurable and, by $(2)$, lies in $L^1(T)$. So, $f\in {\cal A}(T\times E)$, the functional $J_f$ being bounded
in $L^p(T,E)$. By the dominated convergence theorem again, it is likewise clear that
$$\lim_{n\to +\infty}J_{f_n}(u)=J_f(u)$$
for all $u\in L^p(T,E)$. We now show that $J_f$ has no global extrema in $L^p(T,E)$. 
Indeed, fix any $v\in L^p(T,E)$. 
 For each $t\in T$,
set
$$\alpha(t)=\inf_{x\in E}f(t,x)$$
and
$$\beta(t)=\sup_{x\in E}f(t,x).$$
Since $f$ is Carath\'eodory, both $\alpha$ and $\beta$ are measurable functions. Notice that
$$\alpha(t)<f(t,v(t))<\beta(t)$$
for all $t\in D$. Put
$$\gamma_1(t)={{f(t,v(t))-\alpha(t)}\over {2}}$$
and 
$$\gamma_2(t)={{\beta(t)-f(t,v(t))}\over {2}}$$
for all $t\in D$. Further, consider the multifunctions $\Gamma_1, \Gamma_2:D\to 2^E$ defined by
$$\Gamma_1(t)=\{x\in E : f(t,x)\leq \gamma_1(t)\}$$
and
$$\Gamma_2(t)=\{x\in E: f(t,x)\geq\gamma_2(t)\}.$$
By Theorem 6.4 of [3], the multifunctions $\Gamma_1,\Gamma_2$ are measurable, with non-empty closed values. So, by the Kuratowski-Ryll-Nardzewski theorem,
there are two measurable functions $w_1, w_2:D\to E$ such that
$$f(t,w_1(t)\leq \gamma_1(t)$$
and
$$f(t,w_2(t))\geq \gamma_2(t)$$
for all $t\in D$. So, we have
$$f(t,w_1(t))<f(t,v(t))<f(t,w_2(t))$$
for all $t\in D$. Now, fix a measurable subset $A$ of $D$, with $0<\mu(A)<+\infty$, so that the set $w_1(A)\cup w_2(A)$ is bounded and
extend $w_1, w_2$ to the whole of $T$ defining them equal to $v$ outside of $A$. Hence, $w_1, w_2\in L^p(T,E)$.
 Integrating over $T$, since
the measure of $A$ is positive, we get
$$\int_Tf(t,w_1(t))d\mu<\int_Tf(t,v(t))d\mu<\int_Tf(t,w_2(t))d\mu.$$
Hence, as claimed, $J_f$ has no global extrema. Now, fix $V\in {\cal V}(L^p(T,E))$ and $u\in V$. Notice that $u$ is not a local extremum
for the restriction of $J_f$ to $V$. Indeed, otherwise, in view of Corollary 1 of [2], $u$ would be a global extremum for such a restriction. But, in view of Theorem 2 of [7], we would have
$$\inf_{L^p(T,E)}J_f=\inf_VJ_f$$
and
$$\sup_{L^p(T,E)}J_f=\sup_VJ_f$$
and hence $u$ would be a global extremum of $J_f$, against what we have proved. Next, observe that, by Theorem 2 of [6], the set $V$ is an absolute retract and so, by a classical result, it is locally connected.
At this point, the conclusion follows directly by applying Theorem 3.7 of [5] to the restrictions of $J_{f_n}$ and $J_f$ to $V$.\hfill $\bigtriangleup$
\medskip
We now present some applications of Theorem 1.\par
\smallskip
Given two topological spaces $X, Y$, a multifunction $F:X\to 2^Y$ is said to be sequentially lower semicontinuous at the point $\tilde x\in X$
if, for every $\tilde y\in F(\tilde x)$ and for every sequence $\{x_n\}$ in $X$ converging to $\tilde x$, there exists a sequence $\{y_n\}$ in $Y$
converging to $\tilde y$ such that $y_n\in F(x_n)$ for all $n\in {\bf N}$ large enough. We say that $F$ is sequentially lower semicontinuous if it is
so at each point of $X$.\par
\smallskip
Let us also recall that $F$ is said to be lower semicontinuous at $\tilde x$ if, for every open set $\Omega\subset Y$ meeting $F(\tilde x)$, there is a neighbourhood $U$ of $\tilde x$ such that $F(x)\cap \Omega\neq \emptyset$ for all $x\in U$. We say that $F$ is lower semicontinuous if it is
so at each point of $X$.\par
\smallskip
The following proposition is immediate.\par
\medskip
PROPOSITION 1. - {\it The following assertions hold:\par
\noindent
$(a)$\hskip 5pt  If $X$ is first-countable and $F$ is sequentially lower semicontinuous, then $F$ is lower semicontinuous.\par
\noindent
$(b)$\hskip 5pt If $F$ is lower semicontinuous and with non-empty values, then $F$ is sequentially lower semicontinuous.\par
\noindent
$(c)$\hskip 5pt Let $F$ be sequentially lower semicontinuous, let $S$ be a topological space and let $G:Y\to 2^S$ be sequentially lower
semicontinuous. Then, the composite multifunction $x\to G(F(x))$ is sequentially lower semicontinuous.}\par
\medskip
From Theorem 1, we get\par
\medskip
THEOREM 2. {\it Let $M\in L^1(T)$, $c>0$, $\Psi\in (L^p(T,E))^*$, $g\in {\cal A}(T\times E)$. Assume that $J_g$ is Lipschitzian in $L^p(T,E)$,
with Lipschitz constant strictly less than $\|\Psi\|_{(L^p(T,E)^*}$. Set
$$\Gamma=\{f\in {\cal A}(T\times E) : |f(t,x)|\leq M(t)+c\|x\|^p\}$$
and
$$\Lambda=\{f\in {\cal A}(T\times E) : |f(t,x)|\leq M(t)\hskip 5pt and\hskip 5pt \mu\left(\{t\in T: f(t,\cdot)\hskip 5pt has\hskip 5pt  no\hskip 5pt global\hskip 5pt extrema\}\right)>0\}.$$
Moreover, let $\Phi:\Gamma\to 2^{L^p(T,E)}$ be the multifunction defined by
$$\Phi(f)=\{u\in L^p(T,E) : \Psi(u)=J_f(u)+J_g(u)\}.$$
Then, $\Phi$ is sequentially lower semicontinuous with respect to the topology of pointwise convergence at each point of $\Lambda$.}
\smallskip
PROOF. Fix $f\in \Lambda$ and a sequence $\{f_n\}$ in $\Gamma$ pointwise converging to $f$. Also, fix $u\in \Phi(f)$ and set
$$V=\{v\in L^p(T,E) : \Psi(v)-J_g(v)=\Psi(u)-J_g(u)\}.$$
Notice that $V\in {\cal V}(L^p(T,E))$. 
Indeed, take a set $C\in {\cal F}$, with $0<\mu(C)<+\infty$ and consider the function 
$h:T\to {\bf R}$ defined by
$$h(t)=\cases{{{\Psi(u)-J_g(u)}\over {\mu(C)}} & if $t\in C$\cr & \cr 0 & if $t\in T\setminus C$.\cr}$$
Then, if we take $\tilde g(t,x)=g(t,x)+h(t)$, we have $\tilde g\in {\cal A}(T\times E)$, $J_{\tilde g}$ differs from $J_g$ by a constant and
$$V=\{v\in L^p(T,E) : \Psi(v)=J_{\tilde g}(v)\}$$
which shows the claim. Since $u\in V$, by Theorem 1, there exists a sequence $\{u_n\}$ in $L^p(T,E)$ converging to $u$ such that
$$ \Psi(u_n)-J_g(u_n)=\Psi(u)-J_g(u)$$
and
$$J_{f_n}(u_n)=J_f(u)$$
 for all $n\in {\bf N}$ large enough. But
$$\Psi(u)=J_f(u)+J_g(u)$$
and so
$$\Psi(u_n)=J_{f_n}(u_n)+J_g(u_n)$$
and we are done.\hfill $\bigtriangleup$\par
\medskip
Here is another application of Theorem 1.\par
\smallskip
Let $[a,b]$ be a compact real interval. We denote by $AC([a,b],E)$ the space of all absolutely continuous functions $u:[a,b]\to E$, equipped
with the norm
$$\|u\|_{AC([a,b],E)}=\|u(a)\|_E+\int_a^b\|u'(t)\|_Edt.$$
\medskip
THEOREM 3. - {\it Let
let $\{f_n\}$ be a sequence of real-valued Carath\'eodory functions on $T\times E$ such that, for each $n\in {\bf N}$, one has
$$|f_n(t,x)|\leq M(t)+c\|x\|$$
for a.e. $t\in T$ and for all $x\in E$, where $M\in L^1(T)$ and $c\geq 0$. Assume that $\{f_n\}$ converges pointwise in $T\times E$ to a function $f$ such that
$$|f(t,x)|\leq M(t)$$
for a.e. $t\in T$ and for all $x\in E$.
 Moreover, assume that the set
$$D=\{t\in T: f(t,\cdot)\hskip 5pt has\hskip 5pt  no\hskip 5pt global\hskip 5pt extrema\}$$
has a positive measure. Let $[a,b]$ be a compact real interval.\par
Then, for each $r_0, r_1\in {\bf R}$ , for each $\eta\in E^*\setminus \{0\}$, for each $u\in AC([a,b],E)$ such that
$\eta(u(a))=r_0$, $\eta(u(b))=r_1$, and for each sequence $\{\lambda_n\}$ in ${\bf R}$ converging to $\int_a^bf(t,u'(t))dt$, there exist a sequence
$\{u_n\}$ in $AC([a,b],E)$ converging to $u$ and $\nu\in {\bf N}$ such that
$$\eta(u_n(a))=r_0,$$
$$\eta(u_n(b))=r_1$$
and
$$\int_a^bf_n(t,u_n'(t))dt=\lambda_n$$
for all $n\geq\nu$.}\par
\smallskip
PROOF. Let us apply Theorem 1 with $T=[a,b]$ ($\mu$ being the Lebesgue measure) and $p=1$. 
Consider the functional $\Psi:L^1(T,E)\to {\bf R}$
defined by
$$\Psi(v)=\int_a^b\eta(v(t))dt.$$
Of course, $\Psi\in (L^1(T,E))^*$. Notice that $\eta\circ u\in AC([a,b])$ and $(\eta\circ u)'=\eta\circ u'$. Hence, we have
$$\Psi(u')=\int_a^b\eta(u'(t))dt=\eta(u(b))-\eta(u(a))=r_1-r_0.$$
Since $\Psi^{-1}(r_1-r_0)\in {\cal V}(L^1([a,b],E))$, there exist a sequence $\{v_n\}$ in $L^1([a,b],E)$
converging to $u'$ and $\nu\in {\bf N}$ such that
$$\int_a^b\eta(v_n(t))dt=r_1-r_0$$
and
$$\int_a^bf_n(t,v_n(t))dt=\lambda_n.$$
Consider the continuous linear operator $P:L^1([a,b],E)\to AC([a,b],E)$ defined by
$$P(v)(t)=\int_a^tv(s)ds$$
for all $v\in L^1([a,b],E)$, $t\in [a,b]$, the integral being that of Bochner.
Now, put
$$u_n=P(v_n)+u(a).$$
So, $u_n\in AC([a,b],E)$, $u_n'=v_n$, $\lim_{n\to \infty}u_n=P(u')+u(a)=u$. Finally, we have 
$$\eta(u_n(a))=\eta(u(a))=r_0,$$
$$\eta(u_n(b))=\eta\left(\int_a^bv_n(t)dt\right)+\eta(u(a))=\int_a^b\eta(v_n(t))dt+\eta(u(a))=r_1-r_0+\eta(u(a))=r_1$$
and we are done.\hfill $\bigtriangleup$\par
\medskip
Finally, we focus an application of Theorem 2 concerning certain differential inclusions.\par
\smallskip
If $A$ is a non-empty subset of ${\bf R}^k$ and $r>0$, we set
$$B(A,r)=\{y\in {\bf R}^k : \hbox {\rm dist}(y,A)<r\},$$
$$\bar B(A,r)=\{y\in {\bf R}^k : \hbox {\rm dist}(y,A)\leq r\}.$$
\medskip
THEOREM 4. - {\it Let $M, \Psi, g, \Lambda$ be as in Theorem 2 and let $y_0\in {\bf R}^k$, $r>0$. Let $h:\bar B(y_0,r)\to \Lambda$ be
a sequentially continuous function with respect to the topology of pointwise convergence and let $G:L^p(T,E)\to 2^{{\bf R}^k}$ be a sequentially lower semicontinuous
multifunction with non-empty values.\par
Then, there exist $a>0$ and a Lipschitzian function $\varphi:[0,a]\to \bar B(y_0,r)$, with $\varphi(0)=y_0$, having the following property: for a.e. $s\in
[0,a]$ and for each $\epsilon>0$, there exists $u\in L^p(T,E)$ such that
$$\varphi'(s)\in B(G(u),\epsilon)$$
and
$$\Psi(u)=\int_Th(\varphi(s))(t,u(t))d\mu+\int_Tg(t,u(t))d\mu.$$}\par
\smallskip
PROOF. Consider the multifunction $F:\bar B(y_0,r)\to 2^{\bf R^k}$ defined by
$$F(y)=G(\Phi(h(y)))$$
for all $y\in \bar B(y_0,r)$, where $\Phi$ is as in Theorem 2.  By Theorem 2, the restriction to $\Lambda$ of the multifunction $\Phi$ is sequentially lower semicontinuous with respect to the topology of pointwise convergence, and hence
also the multifunction $F$ is sequentially lower semicontinuous, in view of Proposition 1 ($(c)$). Then, since $\bar B(y_0,r)$ is first-countable, $F$ is lower semicontinuous ($(a)$). In particular, this implies that the function
$\hbox {\rm dist}(0,F(\cdot))$ is upper semicontinuous and hence $\rho=\sup_{y\in \bar B(y_0,r)}\hbox {\rm dist}(0, F(y))<+\infty.$
Now, consider the multifunction $H:\bar B(y_0,r)\to {\bf R}^k$ defined by
$$H(y)=F(y)\cap B(0,\rho+1)$$
for all $y\in \bar B(y_0,r)$. Since $B(0,\rho+1)$ is open, $H$ is lower semicontinuous. Hence, $y\to \overline {H(y)}$ is a lower semicontinuous
multifunction, with non-empty closed values and bounded range. As a consequence, by Theorem 2 of [1], there exist $a>0$ and
an absolutely continuous $\varphi:[0,a]\to \bar B(y_0,r)$, with $\varphi(0)=y_0$, such that
$$\varphi'(s)\in  \overline {H(\varphi(s))}\subseteq \overline {F(\varphi(s)}\cap \bar B(0,\rho+1)$$
for a.e. $s\in [0,a]$. Clearly, the function $\varphi$ satisfies the thesis, and we are done.\hfill $\bigtriangleup$.
\medskip
REMARK 1. - We are not aware of known results close enough to Theorem 1 so that a proper comparison can be made.\par
\bigskip
\bigskip
{\bf Acknowledgements:} The author has been supported  by the Gruppo Nazionale per l'Analisi Matematica, la Probabilit\`a e 
le loro Applicazioni (GNAMPA) of the Istituto Nazionale di Alta Matematica (INdAM).
\vfill\eject
\centerline {\bf References}\par
\bigskip
\bigskip
\noindent
[1]\hskip 5pt A. BRESSAN, {\it On differential relations with lower continuous right-hand side. An existence theorem}, J. Differential Equations {\bf 37} (1980), 89-97.\par
\smallskip
\noindent
[2]\hskip 5pt E. GINER, {\it Minima sous contrainte, de fonctionnelles int\'egrales}, C. R. Acad. Sci. Paris S\'er. I Math., {\bf 321} (1995), 429-431.\par
\smallskip
\noindent
[3]\hskip 5pt C. J. HIMMELBERG, {\it Measurable relations},  Fund. Math., {\bf 87} (1975), 53-72.\par
\smallskip
\noindent
[4]\hskip 5pt A. IONESCU TULCEA and C. IONESCU TULCEA, {\it Topics in the theory of lifting}, Springer-Verlag, 1969.\par
\smallskip
\noindent
[5]\hskip 5pt B. RICCERI, {\it On multiselections}, Matematiche, {\bf 38} (1983), 221-235.\par
\smallskip
\noindent
 [6]\hskip 5pt B. RICCERI, {\it Structure, approximation et d\'ependance
continue des solutions de certaines \'equations non lin\'eaires},
C. R. Acad. Sci. Paris, S\'erie I, {\bf 305} (1987), 45-47.\par
\smallskip
\noindent
[7]\hskip 5pt B. RICCERI, {\it Further considerations on a variational property of integral functionals}, J. Optim. Theory Appl., {\bf 106}
(2000), 677-681.\par
\bigskip
\bigskip
\bigskip
\bigskip
Department of Mathematics and Informatics\par
University of Catania\par
Viale A. Doria 6\par
95125 Catania, Italy\par
{\it e-mail address}: ricceri@dmi.unict.it
\bye